\documentclass[12pt]{amsart}

\usepackage{amsmath,xspace,amssymb,mathrsfs}
\usepackage{color}

\input xy
\xyoption{all}
\xyoption{2cell}
\UseAllTwocells
\CompileMatrices

\newcommand{\Spec}{\operatorname{Spec}}
\renewcommand{\phi}{\varphi}

\newcommand{\Ker}{\operatorname{Ker}}

\newcommand{\Min}{\operatorname{Min}}
\newcommand{\Spp}{\operatorname{Spp}}

\newcommand{\Max}{\operatorname{Max}}

\newcommand{\Ann}{\operatorname{Ann}}

\newtheorem{proposition}{Proposition}[section]
\newtheorem{lemma}[proposition]{Lemma} 
\newtheorem{corollary}[proposition]{Corollary}
\newtheorem{theorem}[proposition]{Theorem}

\newtheorem{prop-def}[proposition]{Proposition and definition}

\theoremstyle{definition}
\newtheorem{definition}[proposition]{Definition}
\newtheorem{example}[proposition]{Example}

\newtheorem{remark}[proposition]{Remark}

\begin{document}

\title[N-pure ideals and mid rings]{N-pure ideals and mid rings}

\author[M. Aghajani]{Mohsen Aghajani}
\address{ Department of Mathematics, Faculty of Basic Sciences, University of Maragheh \\
P. O. Box 55136-553, Maragheh, Iran. }

\email{aghajani14@gmail.com}


\begin{abstract}
In this paper, we introduce the concept of N-pure ideal as a generalization of pure ideal. Using this concept, a new and interesting type of rings is presented, we call it mid ring. Also, we provide new characterizations for von Neumann regular and zero-dimensional rings. Moreover, some results about mp-ring are given. Finally, a characterization for mid rings is provided. Then it is shown that the class of mid rings is strictly between the class of reduced mp-rings (p.f. rings) and the class of mp-rings.\\
\end{abstract}

\maketitle

\section{Introduction}

In this  paper, all rings are commutative with identity.
This paper is devoted to study interesting class of ideals which are called N-pure ideals. This notion is a generalization of pure ideal. The theme that encouraged us to study N-pure ideals is that we can use this concept to characterize some different rings. Moreover, this tool lets us to introduce a new type of rings which is completely different with mp-rings and reduced mp-rings. An ideal $I$ of a ring $R$ is said to be a pure ideal if the ring map $R\rightarrow R/I$ is flat, or equivalently for each $a\in I$ there exists $b\in I$ such that $a(1-b)=0$. Pure ideal is an important tool in the study of some areas of ring theory. Pure notion was studied in some works, e.g. \cite{Aghajani and Tarizadeh1}, \cite{Al Ezeh}, \cite{Borceux}, \cite{De Marco} and \cite{Tarizadeh and Aghajani}. \\
In \S 2 we study the notion of N-pure ideal. At first, some basic properties of N-pure ideals are provided. Also, it is shown that the class of N-pure ideals is strictly greater than the class of pure ideals. These two class are the same if and only if the ring is reduced, see Proposition \ref{Lemma IX}. In the following, we provide a characterization for N-pure ideals Theorem \ref{Theorem I}.\\
In \S 3 using the concepts pure and N-pure, we provide characterizations for some well known rings. Recall that a ring $R$ is said to be a mp-ring if each prime ideal of $R$ contains a unique minimal prime ideal of $R$. These rings are studied in \cite{Aghajani and Tarizadeh1} extensively. We describe such ring in terms of N-pure ideals, see Theorem \ref{Theorem II}. Also, a new characterization is given for zero-dimensional rings, see Theorem \ref{Theorem III}. Then, von Neumann regular rings are characterized in terms of pure ideals, see Corollary \ref{Theorem VIII}.\\
Finally, in \S 4 a new ring is introduced which it is called mid ring. We obtain some results about this ring. As an important result of this section, we provide a characterization for mid rings, see \ref{Theorem VII}. Especially, we prove that a ring $R$ is a mid ring if and only if $R_{\mathfrak{p}}$ is a primary ring for every prime ideal $\mathfrak{p}$ of $R$. Also, a new characterizations for p.p. rings is given, see Theorem \ref{Theorem IX} and Corollary \ref{Corollary V}. In the following, the N-pure prime ideals of mid rings are identified, see Theorem \ref{Theorem IV}. Finally, we prove that the class of mid rings is strictly between the class of reduced mp-rings and the class of mp-rings, see Remark \ref{Remark I}, Theorem \ref{Theorem X} and Example \ref{Example I}.\\

Now we recall some notions which use in this paper. By p.p. ring we mean a ring whose each of its principal ideal is projective. A ring $R$ is said to be a p.f. ring if for each $a\in R$, $\Ann(a)$ is a pure ideal. Clearly every p.p. ring is a p.f. ring. A ring $R$ is a Gpf-ring if for every $a\in R$
there exists $n\geqslant1$ such that $\Ann(a^{n})$ is a pure ideal. We denote nil-radical (Jacobson radical) of a ring $R$ by $\mathfrak{N}$ (resp. $\mathfrak{J}$). Also, nil-radical (Jacobson radical) of ring $R_{\mathfrak{p}}$ is denoted by $\mathfrak{N}(R_{\mathfrak{p}})$ (resp. $\mathfrak{J}(R_{\mathfrak{p}})$), for a prime ideal $\mathfrak{p}$ of $R$. \\

\section{N-pure ideals}

We begin this section with the following definition.\\

\begin{definition}\label{Def I}
An ideal $I$ of a ring $R$ is called \emph{N-pure} if for every $a\in I$ there exists $b\in I$ such that $a(1-b) \in \mathfrak{N}$.\\
\end{definition}

\begin{remark}\label{Remark IV}
 It is clear that every pure ideal is N-pure. Moreover, if $I$ is a pure ideal of  ring $R$, then its radical is a N-pure ideal. In particular, the nil-radical of $R$ is a N-pure ideal. Therefore, if $R$ is a non-reduced ring, it is straightforward to check that $\mathfrak{N}$ is a N-pure ideal which is not pure. The following result describes reduced rings in terms of pure and N-pure ideals.\\
\end{remark}

\begin{proposition}\label{Lemma IX}Let $R$ be a ring. Then $R$ is a reduced ring if and only if every N-pure ideal is a pure ideal.\\
\end{proposition}

{\bf Proof.} Let $I$ be a N-pure ideal and $a \in I$. Then there exists $b \in I$ and such that $a(1-b) \in \mathfrak{N}$. Thus $a(1-b)=0$. Therefore, $I$ is a pure ideal. The converse follows from Remark \ref{Remark IV}. $\Box$\\

\begin{proposition}\label{Proposition I} Let $I$ be an ideal of a ring $R$. Then $I$ is a N-pure ideal if and only if for each $a \in I$ there exist $n\geqslant1$ and $b \in I$ such that $a^{n}(1-b)=0$.\\
\end{proposition}

{\bf Proof.} Let $I$ be a N-pure ideal and $a \in I$. Then there exists $c \in I$ such that $a(1-c) \in \mathfrak{N}$. Thus there exists $n\geqslant1$ such that $(a(1-c))^{n}=0$ and so $a^{n}(1-b)=0$ for some $b \in I$. Conversely, if  $a^{n}(1-b)=0$, then we have  $a(1-b) \in \mathfrak{N}$ and so $I$ is N-pure. $\Box$ \\

\begin{lemma}\label{Lemma III} Let $I$ be an ideal of a ring $R$. Then $I$ is a N-pure ideal if and only if $(I+\mathfrak{N})/\mathfrak{N}$ is a pure ideal.\\
\end{lemma}

{\bf Proof.} Assume that $I$ be a N-pure ideal. If $a \in I$ then there exists $b \in I$ such that $a(1-b) \in \mathfrak{N}$. Thus $(I+\mathfrak{N})/\mathfrak{N}$ is a pure ideal. The converse case is straightforward.  $\Box$ \\

Recall that a subset E of $\Spec(R)$ is said to be stable under the generalization if for
any two prime ideals $\mathfrak{p}$ and $\mathfrak{q}$ of $R$ with $\mathfrak{p} \subset \mathfrak{q}$, if $\mathfrak{q} \in E$ then $\mathfrak{p} \in E$. The following result provides a characterization for N-pure ideals. \\

\begin{theorem}\label{Theorem I} Let $I$ be an ideal of a ring $R$. Then the following conditions are equivalent. \\
$\mathbf{(i)}$ $I$ is a N-pure ideal. \\
$\mathbf{(ii)}$ For every $a_{1},...,a_{n} \in I$ there exist $b \in I$ and $t\geqslant1$ such that $a_{k}^{t}=a_{k}^{t}b$ for all $k=1,...,n$.\\
$\mathbf{(iii)}$ For every $a \in I$ there exists $t\geqslant1$ such that $\Ann(a^{t})+I=R$. \\
$\mathbf{(iv)}$ $\sqrt{I}=\{a\in R|\exists n \geqslant1, \Ann(a^{n})+I=R\}$.\\
$\mathbf{(v)}$ $\sqrt{I}$ is a N-pure ideal.\\
$\mathbf{(vi)}$ There exists a unique pure ideal $J$ such that $\sqrt{I}=\sqrt{J}$.\\
\end{theorem}

{\bf Proof.} $\mathbf{(i)}\Rightarrow\mathbf{(ii)}:$ Since $I$ is N-pure, then there exist $b_{k} \in I$ and $t_{k}\geqslant1$ such that $a_{k}^{t_{k}}=a_{k}^{t_{k}}b_{k}$. Setting $t:=\max\{t_{k}| 1\leq k \leq n\}$ and let $b \in I$ where $1-b=\prod\limits_{k=1}^{n}(1-b_{k})$. Clearly, we have $b \in I$ and $a_{k}^{t}=a_{k}^{t}b$ and so the assertion is proved.\\
$\mathbf{(ii)}\Rightarrow\mathbf{(iii)}:$ Let $a \in I$. Then by assumption there exist $b \in I$ and $t\geqslant1$ such that $a^{t}=a^{t}b$. Thus we have $1-b \in \Ann(a^{t})$ and so $\Ann(a^{t})+I=R$.\\
$\mathbf{(iii)}\Rightarrow\mathbf{(iv)}$ Let $a \in \sqrt{I}$. Then there exists $m\geqslant1$ such that $a^{m} \in I$ and so by hypothesis, we have $\Ann(a^{t})+I=R$ for some $t\geqslant1$. Conversely, if $\Ann(a^{n})+I=R$ for some $n\geqslant1$, then we have $a^{n} \in I$ and so $a \in \sqrt{I}$.\\
$\mathbf{(iv)}\Rightarrow\mathbf{(v)}:$ Let $a\in \sqrt{I}$. Then there exists $n\geqslant1$ such that $\Ann(a^{n})+I=R$. Thus we have $c+d=1$ for some $c\in \Ann(a^{n})$ and $d\in I$. Hence $a^{n}(1-d)=0$ and so $\sqrt{I}$ is a N-pure ideal.\\
$\mathbf{(v)}\Rightarrow\mathbf{(vi)}:$ Let $\mathfrak{p} \in V(I)$ and $\mathfrak{q} \subset \mathfrak{p}$. If $\mathfrak{q} \notin V(I)$, then there exists $a \in I\setminus \mathfrak{q}$. Since $\sqrt{I}$ is N-pure, then there are $b \in \sqrt{I}$ and $n\geqslant1$ such that $a^{n}(1-b)=0$. Hence, $1-b \in \mathfrak{p}$ and so $1 \in \mathfrak{p}$ which is a contradiction. Therefore, $V(I)$ is stable under the generalization. Now by \cite[Theorem 3.2]{Ta}, there exists a unique pure ideal $J$ such that $\sqrt{I}=\sqrt{J}$.\\
$\mathbf{(vi)}\Rightarrow\mathbf{(i)}:$ Let $a \in I$. Then there exists $t\geqslant1$ such that $a^{t}\in J$. Thus there exists $b \in J$ such that $a^{t}(1-b)=0$, since $J$ is a pure ideal. On the other hand, we have $b^{s}\in I$ for some $s\geqslant1$. Therefore, $a^{t}(1-b^{s})=0$ and so $I$ is a N-pure ideal. $\Box$\\

\begin{corollary}\label{Corollary III} Let $I$ be an ideal of a ring $R$. Then $I$ is a N-pure ideal if and only if $I^{n}$ is a N-pure ideal for all $n\geqslant1$. $\Box$ \\
\end{corollary}

\begin{theorem}\label{Theorem VI} Let $(I_{k})$ be a family of N-pure ideals of a ring $R$. Then the following statements hold.\\
$\mathbf{(i)}$ $\sum\limits_{k}I_{k}$ is a N-pure ideal.\\
$\mathbf{(ii)}$ If $I$ and $J$ are N-pure ideals of $R$, then $IJ$ and $I \cap J$ are N-pure ideals.\\
\end{theorem}

{\bf Proof.} $\mathbf{(i)}$ Let $a= \sum\limits_{j=1}^{t}a_{j}\in \sum\limits_{k}I_{k}$ where $a_{j} \in I_{j}=I_{k_{j}}$. Then there exist $b_{j} \in I_{j}$ and $n_{j}\geqslant1$ such that $a_{j}^{n_{j}}(1-b_{j})=0$. Setting $1-b=\prod\limits_{j=1}^{t}(1-b_{j})$ and $n=\sum\limits_{j=1}^{t}n_{j}$. It is easy to see that $b \in \sum\limits_{k}I_{k}$. But we have $a^{n}(1-b)=\big(\sum\limits_{j=1}^{t}a_{j}\big)^{n}(1-b)=0$. Thus $\sum\limits_{k}I_{k}$ is a N-pure ideal. $\mathbf{(ii)}$ Since the product of two pure ideals is a pure ideal, then the assertions follow from Theorem \ref{Theorem I}(vi). $\Box$ \\

\begin{lemma}\label{Lemma VI} Let $R$ be a ring and $I$ be an ideal of $R$. If for every $x \in R$ there exists $n\geqslant1$ such that $Rx^{n}\cap I=x^{n}I$, then $I$ is a N-pure ideal.\\
\end{lemma}

{\bf Proof.} Let $a \in I$. Then there exists $n\geqslant1$ such that $Ra^{n}\bigcap I=a^{n}I$. So there exists $b\in I$ such that $a^{n}=a^{n}b$. Therefore, $I$ is a N-pure ideal. $\Box$ \\

\section{some characterizations of rings}

We begin this section with the following remark.\\

\begin{remark}\label{Remark III}Let $R$ be a ring. If $\mathfrak{p}$ is a minimal prime ideal of $R$, then $\sqrt{\Ker\pi_{\mathfrak{p}}}=\mathfrak{p}$, because if $\Ker\pi_{\mathfrak{p}}\subseteq \mathfrak{q}$ and $\mathfrak{p}\nsubseteq\mathfrak{q}$, then there exists $a\in \mathfrak{p}\setminus \mathfrak{q}$. So there exists $s\in R\setminus \mathfrak{p}$ such that $sa\in N(R)$. Hence $a^{n}\in \Ker\pi_{\mathfrak{p}}$ for some $n\geqslant1$ and so $a^{n}\in \mathfrak{q}$ which is a contradiction. Therefore we have $V(\Ker\pi_{\mathfrak{p}})=V(\mathfrak{p})$.
\end{remark}
Recall that a ring $R$ is said to be a \emph{$NJ$-ring} if  its nil-radical and jacobson radical are the same ($\mathfrak{N}=\mathfrak{J}$). Also a ring $R$ is called a \emph{semiprimitive} ring if $\mathfrak{J}=0$. The following result provides characterizations for such rings.\\

\begin{proposition}\label{proposition III} The following statements hold.\\
$\mathbf{(i)}$ $R$ is a $NJ$-ring if and only if $\mathfrak{J}$ is a N-pure ideal.\\
$\mathbf{(ii)}$ $R$ is a semiprimitive ring if and only if $\mathfrak{J}$ is a pure ideal. \\
\end{proposition}

{\bf Proof.} $\mathbf{(i)}$ By Remark \ref{Remark IV} $\mathfrak{J}$ is a N-pure ideal of $R$. Conversely, if $a \in \mathfrak{J}$, then there exist $n\geqslant1$ and $b \in \mathfrak{J}$ such that $a^{n}=a^{n} b$. Hence we have $a^{n}(1-b)=0$. Thus $a^{n}=0$ and so $a \in \mathfrak{N}$. Then $R$ is a $NJ$-ring. $\mathbf{(ii)}$ It is clear. Let $a\in \mathfrak{J}$. Then there exists $b\in \mathfrak{J}$ such that $a(1-b)=0$. Hence $a=0$ and so $\mathfrak{J}=0$. Therefore $R$ is a semiprimitive ring. $\Box$ \\

In the following result, we provide a characterization for mp-rings. \\

\begin{theorem}\label{Theorem II} Let $R$ be a ring. Then the following statements are equivalent.\\
$\mathbf{(i)}$ $R$ is a $mp$-ring.\\
$\mathbf{(ii)}$ If $ab=0$, then there exists $n\geqslant1$ such that $\Ann(a^{n})+\Ann(b^{n})=R$.\\
$\mathbf{(iii)}$ Every minimal prime ideal of $R$ is a N-pure ideal.\\
$\mathbf{(iv)}$ For every minimal prime ideal $\mathfrak{p}$ of $R$, $\Ker\pi_{\mathfrak{p}}$ is a N-pure ideal.\\
$\mathbf{(v)}$ For every prime ideal $\mathfrak{p}$ of $R$, $\Ker\pi_{\mathfrak{p}}$ is a N-pure ideal.\\
\end{theorem}

{\bf Proof.} For $\mathbf{(i)}\Leftrightarrow\mathbf{(ii)}$ See \cite[Theorem 6.2]{Aghajani and Tarizadeh1}.\\
$\mathbf{(ii)}\Rightarrow\mathbf{(iii)} :$ Let $\mathfrak{p}\in \Min(R)$. If $a\in \mathfrak{p}$, then there exists $b\in R\setminus\mathfrak{p}$ such that $ab\in N(R)$ and so by hypothesis, there exists $n\geqslant1$ such that $\Ann(a^{n})+\Ann(b^{n})=R$. Then there exist $c\in \Ann(a^{n})$ and $d\in \Ann(b^{n})$ such that $c+d=1$. Thus we have $a^{n}(1-d)=0$ and so $\mathfrak{p}$ is a N-pure ideal, since $d\in \mathfrak{p}$.\\
$\mathbf{(iii)}\Leftrightarrow\mathbf{(iv)} :$ It follows from Theorem \ref{Theorem I} and Remark \ref{Remark III}.\\
$\mathbf{(iv)}\Rightarrow\mathbf{(i)}$ Let $\mathfrak{p}$ and $\mathfrak{q}$ be distinct minimal prime ideals of $R$. Then by Lemma \cite[Lemma 3.2]{Aghajani and Tarizadeh1}, there exist $a\in R\setminus \mathfrak{p}$ and $b\in R\setminus \mathfrak{q}$ such that $ab=0$. Hence $b\in \Ker\pi_{\mathfrak{p}}$.  Thus there exist $n\geqslant1$ and $c\in \Ker\pi_{\mathfrak{p}}$ such that $b^{n}(1-c)=0$, since $\Ker\pi_{\mathfrak{p}}$ is a N-pure ideal. So $1-c\in \mathfrak{q}$. Therefore, $\mathfrak{p}+\mathfrak{q}=R$ and so $R$ is a mp-ring by \cite[Theorem 6.2(ii)]{Aghajani and Tarizadeh1}.\\
$\mathbf{(ii)}\Rightarrow\mathbf{(v)} :$ Let $ab\in \sqrt{\Ker\pi_{\mathfrak{p}}}$ and $a\notin \sqrt{\Ker\pi_{\mathfrak{p}}}$. Then there exists $n\geqslant1$ such that $a^{n}b^{n}\in \Ker\pi_{\mathfrak{p}}$. Hence there exists $s\in R\setminus \mathfrak{p}$ such that $sa^{n}b^{n}=0$. Therefore by hypothesis, there exists $m\geqslant1$ such that $\Ann(a^{mn})+\Ann(s^{m}b^{mn})=R$. Thus $\Ann(s^{m}b^{mn})\nsubseteq \mathfrak{p}$ and so there exists $t\in R\setminus \mathfrak{p}$ such that $ts^{m}b^{mn}=0$. Then $b^{mn}\in \Ker\pi_{\mathfrak{p}}$ and so $b\in \sqrt{\Ker\pi_{\mathfrak{p}}}$. Therefore $\sqrt{\Ker\pi_{\mathfrak{p}}}$ is a minimal prime ideal of $R$. Now, the assertion follows from Theorem \ref{Theorem I} and $(ii)\Rightarrow(iii)$.\\
$\mathbf{(v)}\Rightarrow\mathbf{(i)} :$ Let $\mathfrak{p}$ and $\mathfrak{q}$ be distinct minimal prime ideals of $R$. Then there exist $a\in R\setminus \mathfrak{p}$ and  $b\in R\setminus \mathfrak{q}$ such that $ab=0$. Thus $b\in \Ker\pi_{\mathfrak{p}}$. So by assumption, there exist $n\geqslant1$ and $c\in \Ker\pi_{\mathfrak{p}}$ such that $b^{n}(1-c)=0$. Hence $1-c \in \mathfrak{q}$ and so $\mathfrak{p}+\mathfrak{q}=R$. This means that $R$ is a mp-ring. $\Box$ \\

It is well known that a ring $R$ is von Neumann regular if and only if it is a reduced zero-dimensional ring. Clearly, a ring $R$ is a semiprimitive local ring if and only if it is a field. In the following result using N-pure concept, a new characterization for zero-dimensional rings is given. \\

\begin{theorem}\label{Theorem III} Let $R$ be a ring. Then the following statements are equivalent.\\
$\mathbf{(i)}$ $R$ is a zero-dimensional ring.\\
$\mathbf{(ii)}$ Every ideal of $R$ is a N-pure ideal.\\
$\mathbf{(iii)}$ Every principal ideal of $R$ is a N-pure ideal.\\
$\mathbf{(iv)}$ Every maximal ideal of $R$ is a N-pure ideal.\\
$\mathbf{(v)}$ $\sqrt{\Ker\pi_{\mathfrak{m}}}=\mathfrak{m}$
for all $\mathfrak{m}\in \Max(R)$.\\
$\mathbf{(vi)}$ $R_{\mathfrak{p}}$ is a $NJ$-ring for all  $\mathfrak{p}\in \Spec(R)$.\\
$\mathbf{(vii)}$ $R_{\mathfrak{m}}$ is a $NJ$-ring for all  $\mathfrak{m}\in \Max(R)$.\\
\end{theorem}

{\bf Proof.} $\mathbf{(i)}\Rightarrow\mathbf{(ii)} :$ We know that $R/\mathfrak{N}$ is a von Neumann regular ring. Then if $I$ is an ideal of $R$ and $a\in I$, then there exists $b\in R$ such that $(a-a^{2}b)\in \mathfrak{N}$. Thus $I$ is a N-pure ideal.\\
$\mathbf{(ii)}\Rightarrow\mathbf{(iii)} :$ There is nothing to prove.\\
$\mathbf{(iii)}\Rightarrow\mathbf{(iv)} :$ It is obvious.\\
$\mathbf{(iv)}\Rightarrow\mathbf{(v)} :$ Let $a\in\mathfrak{m}$. Then there exist $n\geqslant1$ and $b\in\mathfrak{m}$ such that $a^{n}(1-b)=0$. Thus $a^{n}\in \Ker\pi_{\mathfrak{m}}$ and so $\sqrt{\Ker\pi_{\mathfrak{m}}}=\mathfrak{m}$.\\
$\mathbf{(v)}\Rightarrow\mathbf{(vi)} :$ Let $a/s\in \mathfrak{J}(R_{\mathfrak{p}})$. Assume that $\mathfrak{m}$ is a maximal ideal which $\mathfrak{p}\subseteq\mathfrak{m}$. Thus $a\in \mathfrak{m}$ and so there exist a natural number $n\geqslant1$ and $t\in R\setminus\mathfrak{m}$ such that $ta^{n}=0$. Hence $a^{n}/s^{n}=0$ in $R_{\mathfrak{p}}$. Therefore we have $\mathfrak{J}(R_{\mathfrak{p}})=\mathfrak{N}(R_{\mathfrak{p}})$.\\
$\mathbf{(vi)}\Rightarrow\mathbf{(vii)} :$ It is clear.\\
$\mathbf{(vii)}\Rightarrow\mathbf{(i)} :$ Let $\mathfrak{p}$ and $\mathfrak{q}$ be prime ideals of $R$ which $\mathfrak{p}\subseteq\mathfrak{q}$. If $\mathfrak{m}$ is a maximal ideal of $R$ containing $\mathfrak{q}$, then we have $\mathfrak{N}(R_{\mathfrak{m}})\subseteq\mathfrak{p}R_{\mathfrak{m}}
\subseteq\mathfrak{q}R_{\mathfrak{m}}
\subseteq\mathfrak{m}R_{\mathfrak{m}}=\mathfrak{J}(R_{\mathfrak{m}})$. Hence $\mathfrak{p}R_{\mathfrak{m}}=\mathfrak{q}R_{\mathfrak{m}}$ and so $R$ is a zero-dimensional ring. $\Box$\\

 In the following result, using the previous theorem and pure notion, a characterization for von Neumann regular rings is given which some of its conditions are well known.\\

\begin{corollary}\label{Theorem VIII}Let $R$ be a ring. Then the following statements are equivalent.\\
$\mathbf{(i)}$ $R$ is a von Neumann regular ring.\\
$\mathbf{(ii)}$ Every ideal of $R$ is a pure ideal.\\
$\mathbf{(iii)}$ Every principal ideal of $R$ is a pure ideal.\\
$\mathbf{(iv)}$ Every maximal ideal of $R$ is a pure ideal.\\
$\mathbf{(v)}$ $\Ker\pi_{\mathfrak{m}}=\mathfrak{m}$ for all $\mathfrak{m}\in \Max(R)$.\\
$\mathbf{(vi)}$ $R_{\mathfrak{p}}$ is a semiprimitive ring for all  $\mathfrak{p}\in \Spec(R)$.\\
$\mathbf{(vii)}$ $R_{\mathfrak{m}}$ is a semiprimitive ring for all  $\mathfrak{m}\in \Max(R)$.\\
\end{corollary}

{\bf Proof.} $\mathbf{(i)}\Rightarrow\mathbf{(ii)} :$  This follows from Proposition \ref{Lemma IX} and Theorem \ref{Theorem III}.\\
$\mathbf{(ii)}\Rightarrow\mathbf{(iii)} :$ There is nothing to prove.\\ $\mathbf{(iii)}\Rightarrow\mathbf{(iv)} :$ It is obvious.\\ $\mathbf{(iv)}\Rightarrow\mathbf{(v)} :$ Suppose $\mathfrak{m}\in\Max(R)$ and $a\in \mathfrak{m}$. Then by the hypothesis, there exists $b\in \mathfrak{m}$ such that $a(1-b)=0$. Thus $a\in \Ker\pi_{\mathfrak{m}}$ and so $\Ker\pi_{\mathfrak{m}}=\mathfrak{m}$.\\
$\mathbf{(v)}\Rightarrow\mathbf{(vi)} :$ It is easy to see that $R$ is reduced. Then by Theorem \ref{Theorem III} $R_{\mathfrak{p}}$ is a semiprimitive ring.\\ $\mathbf{(vi)}\Rightarrow\mathbf{(vii)} :$ There is nothing to prove.\\ $\mathbf{(vii)}\Rightarrow\mathbf{(i)} :$ By Theorem \ref{Theorem III}, it suffices to show that $R$ is reduced.  Let $a\in \mathfrak{N}$. Then $a\in\mathfrak{m}$ for all $\mathfrak{m}\in \Max(R)$. Hence by the hypothesis, there exists $b_{\mathfrak{m}}\in R\setminus\mathfrak{m}$ such that $b_{\mathfrak{m}}a=0$. Thus $I=\big(b_{\mathfrak{m}}: \mathfrak{m}\in\Max(R)\big) $ is equal to $R$. Therefore, we have $1=\sum\limits_{i=1}^{n}r_{i}b_{i}$ where $b_{i}=b_{\mathfrak{m_{i}}}$ and $r_{i}\in R$. Hence $a=\sum\limits_{i=1}^{n}r_{i}b_{i}a=0$ and so $R$ is a reduced ring. $\Box$\\

Recall that a proper ideal $P$ of a ring $R$ is said to be purely-prime if for pure ideals $I$ and $J$ of $R$ with $IJ\subseteq P$ we have either $I\subseteq P$ or $J\subseteq P$ \cite{Tarizadeh and Aghajani}. The set of all purely-prime ideals of a ring $R$ is called the pure spectrum of $R$ and is denoted by $\Spp(R)$. \\

\begin{corollary}\label{Corollary II} Let $R$ be a ring. Then $R$ is a von Neumann regular ring if and only if the prime spectrum of $R$ and the pure spectrum of it are the same.\\
\end{corollary}

{\bf Proof.} This is an immediate consequence of \ref{Theorem VIII}. $\Box$\\

\section{Mid rings}

In this section, we introduce a new class of rings and study some basic properties of it.\\

\begin{definition}\label{Def II}
A ring $R$ is called a \emph{mid ring} if for every $a \in R$, $\Ann(a)$ is a N-pure ideal.
\end{definition}

\begin{proposition}\label{Proposition VI} Let $R$ be a mid ring. Then $R_{\mathfrak{p}}$ is a mid ring for all $\mathfrak{p} \in \Spec(R)$.\\
\end{proposition}

{\bf Proof.} Let $a/s \in R_{\mathfrak{p}}$ and $b/t \in \Ann(a/s)$, then there exists $u \in R\backslash \mathfrak{p}$ such that $uba=0$. Thus there exist $c \in \Ann(a)$ and $ m\geqslant1$ such that $(ub)^{m}(1-c)=0$ and so we have $(ub)^{m}/(ut)^{m}(1-c)=0$. Therefore, $\Ann(a/s)$ is N-pure and the claim is proved. $\Box$ \\

\begin{lemma}\label{Lemma VIII} Let $R$ be a mid ring. If $I$ is a pure ideal, Then $R/I$ is a mid ring.\\
\end{lemma}

{\bf Proof.} Let $a \in R\backslash I$. Then $\Ann(a)$ is a N-pure ideal. If $b+ I \in \Ann(a+I)$, then we have $ab \in I$ and so there exists $c \in I$ such that $ab(1-c)=0$. Thus there exist $d \in \Ann(a)$ and $m\geqslant1$ such that $b^{m}(1-c)^{m}(1-d)=0$. Therefore $b^{m}(1-d) \in I$ and so $R/I$ is a mid ring. $\Box$ \\

\begin{remark}\label{Remark II}It is straightforward to see that if $I_{k}$ is an ideal of ring $R_{k}$, $R=\prod\limits_{k=1}^{n} R_{k}$ and $I=\prod\limits_{k=1}^{n} I_{k}$, then $I$ is a N-pure ideal if and only if every $I_{k}$ is a N-pure ideal. On the other hand, for $a=(a_{k}) \in R$ we have $\Ann(a)=\prod\limits_{k=1}^{n} \Ann(a_{k})$. Then we can obtain the following result.\\
\end{remark}

\begin{proposition}\label{Proposition VII} Let $R=\prod\limits_{k}R_{k}$. If $R$ is a mid ring, then every $R_{k}$ is mid ring. If index set is finite, then the converse holds.\\
\end{proposition}

{\bf Proof.} Let $a_{k_{0}}\in R_{k_{0}}$. Setting  $a=(a_{k})$ where $a_{k} =a_{k_{0}}$ if $k=k_{0}$ and $a_{k} =0$ if $k\neq k_{0}$. Then there exist $b=(b_{k})\in \Ann(a)$ and $n\geqslant1$ such that $a^{n}(1-b)=0$. Thus $a_{k_{0}}^{n}(1-b_{k_{0}})=0$ where $b_{k_{0}}\in \Ann(a_{k_{0}})$.  Then $\Ann(a_{k_{0}})$ is a N-pure ideal of $R_{k_{0}}$ and so $R_{k_{0}}$ is a mid ring. The last assertion follows easily from Remark \ref{Remark II}. $\Box$ \\

\begin{lemma}\label{Lemma II} Let $R$ be a ring. Then  $\Ker\pi_{\mathfrak{p}}$ is a primary ideal for each $\mathfrak{p}\in \Min(R)$.\\
\end{lemma}

{\bf Proof.} Let $\mathfrak{p}$  be a minimal prime ideal of $R$. If $ab \in \Ker\pi_{\mathfrak{p}}$ and $a \notin \Ker\pi_{\mathfrak{p}}$, then there exists $c\in R\setminus \mathfrak{p}$ such that $abc=0$. Since $\Ann(a)\subseteq \mathfrak{p}$, then $bc \in \mathfrak{p}$ and so $b \in \mathfrak{p}$. Thus $\Ker\pi_{\mathfrak{p}}$ is a primary ideal by Remark \ref{Remark III}. $\Box$ \\

Recall that a ring is a primary ring if its zero ideal is a primary ideal. The following result provides a characterization for mid rings.\\

\begin{theorem}\label{Theorem VII} Let $R$ be a ring. Then the following are equivalent.\\
$\mathbf{(i)}$ $R$ is a mid ring.\\
$\mathbf{(ii)}$ If $ab=0$, then there exists $n\geqslant1$ such that $\Ann(a) + \Ann(b^{n})=R$.\\
$\mathbf{(iii)}$ $R_{\mathfrak{p}}$ is a primary ring for all $\mathfrak{p}\in \Spec(R)$.\\
$\mathbf{(iv)}$ $R_{\mathfrak{m}}$ is a primary ring for all $\mathfrak{m}\in \Max(R)$.\\
$\mathbf{(v)}$ $\Ker\pi_{\mathfrak{p}}$ is a pure ideal for all $\mathfrak{p}\in \Spec(R)$.\\
$\mathbf{(vi)}$ $\Ker\pi_{\mathfrak{p}}$ is a pure ideal for all $\mathfrak{p}\in \Min(R)$.\\
$\mathbf{(vii)}$ $\Ker\pi_{\mathfrak{p}}=\Ker\pi_{\mathfrak{q}}$ for prime ideals $\mathfrak{p}$ and $\mathfrak{q}$ with $\mathfrak{p}\subseteq\mathfrak{q}$.\\
$\mathbf{(viii)}$ $\Ker\pi_{\mathfrak{p}}$ is a primary ideal for all $\mathfrak{p}\in \Spec(R)$.\\
$\mathbf{(ix)}$ $\Ker\pi_{\mathfrak{m}}$ is a primary ideal for all $\mathfrak{m}\in \Max(R)$.\\
 
\end{theorem}

{\bf Proof.} $\mathbf{(i)} \Rightarrow \mathbf{(ii)}$ : Let $ab=0$. Then by hypothesis, $\Ann(a)$ is a N-pure ideal and so there are $n\geqslant1$ and $c \in \Ann(a)$ such that $b(1-c) \in \mathfrak{N}$. Thus we have $\Ann(a) + \Ann(b^{n})=R$. \\
$\mathbf{(ii)} \Rightarrow \mathbf{(iii)}$ : Let $(a/s)(b/s')=0$ and $a/s\neq0$ in $R_{\mathfrak{p}}$. Then there exists $t\in R \setminus \mathfrak{p}$ such that $tab=0$. Thus there exists $n\geqslant1$ such that $\Ann(a) + \Ann((tb)^{n})=R$. Therefore $\Ann((tb)^{n})\nsubseteq\mathfrak{p}$ and so $b/s'$ is a nilpotent in $R_{\mathfrak{p}}$.\\
$\mathbf{(iii)} \Rightarrow \mathbf{(iv)}$ : There is nothing to prove.\\
$\mathbf{(iv)} \Rightarrow \mathbf{(ii)}$ : Let $ab=0$ for some $a,b\in R$. Suppose that $\Ann(a) + \Ann(b^{n})\neq R$ for all
$n\geqslant1$. Setting $I:=\sum\limits_{n\geqslant1}\big(\Ann(a)+\Ann(b^{n})\big)$, then there exists a maximal ideal $\mathfrak{m}$ such that $I\subseteq\mathfrak{m}$. Hence we have $a/1\neq0$ and $b/1\nsubseteq \mathfrak{N}(R_{\mathfrak{m}})$ which is a contradiction.\\
$\mathbf{(ii)} \Rightarrow \mathbf{(i)}$ : Let $a\in R$ and $b\in \Ann(a)$. Then there exists $n\geqslant1$ such that $\Ann(a) + \Ann(b^{n})=R$. Hence there exists $c\in \Ann(a)$ such that $b(1-c)\in \mathfrak{N}$ and so $R$ is a mid ring.\\
$\mathbf{(ii)}\Rightarrow\mathbf{(v)} :$ Let $\mathfrak{p}$ a prime ideal of $R$ and $a\in \Ker\pi_{\mathfrak{p}}$. Then there exists $b\in R\setminus\mathfrak{p}$ such that $ab=0$. Thus by hypothesis there is $n\geqslant1$ such that $\Ann(a) + \Ann(b^{n})=R$. Therefore we have $c+d=1$ for some $c\in \Ann(a)$ and $d\in \Ann(b^{n})$. Hence $d\in \Ker\pi_{\mathfrak{p}}$ and $a(1-d)=0$. Then $\Ker\pi_{\mathfrak{p}}$ is a pure ideal of $R$. \\
$\mathbf{(v)}\Rightarrow\mathbf{(vi)} :$ There is nothing to prove.\\
$\mathbf{(vi)}\Rightarrow\mathbf{(vii)} :$ Let $\mathfrak{m}$ be a maximal ideal of $R$ containing minimal prime ideal $\mathfrak{p}$. Then we have $\Ker\pi_{\mathfrak{m}} \subseteq \Ker\pi_{\mathfrak{p}}$. Let $a \in \Ker\pi_{\mathfrak{p}}$. Since $\Ker\pi_{\mathfrak{p}}$ is a pure ideal, there exists $b \in \Ker\pi_{\mathfrak{p}}$ such that $a(1-b)=0$. Thus $1-b \notin \mathfrak{m}$ and so $a \in \Ker\pi_{\mathfrak{m}}$. This yields that $\Ker\pi_{\mathfrak{m}}=\Ker\pi_{\mathfrak{p}}$. \\
$\mathbf{(vii)}\Rightarrow\mathbf{(viii)} :$ It follows from Lemma \ref{Lemma II}.\\
$\mathbf{(viii)} \Rightarrow \mathbf{(ix)}$ : It is obvious.\\
$\mathbf{(ix)} \Rightarrow \mathbf{(ii)}$ : Assume that $ab=0$. we claim that $\Ann(a) + \Ann(b^{m})=R$ for some $m\geqslant1$. Otherwise, there exists $\mathfrak{m}\in \Max(R)$ such that $\sum\limits_{n\geqslant1}\big(\Ann(a) + \Ann(b^{n})\big)\subseteq \mathfrak{m}$. Thus $a\notin \Ker\pi_{\mathfrak{m}}$ and $b\notin \sqrt{\Ker\pi_{\mathfrak{m}}}$ which is a contradiction.
$\Box$ \\

\begin{remark}\label{Remark V} It is easy to see that if $R$ is a primary ring, then $R_{\mathfrak{p}}$ is a primary ring for all $\mathfrak{p}\in \Spec(R)$. Then by Theorem \ref{Theorem VII}, every primary ring is a mid ring.\\
Now let $R$ be a $Gpf$-ring. Thus for each $a\in R$ there exists $n\geqslant1$, whenever $a^{n}b=0$, then $\Ann(a^{n}) + \Ann(b)=R$. Because if $a\in R$, then there exists $n\geqslant1$ such that $\Ann(a^{n})$ is a pure ideal. Thus easily we have $\Ann(a^{n}) + \Ann(b)=R$ for each $b\in \Ann(a^{n})$. \\
\end{remark}

Now we can obtain the next result.\\

\begin{lemma}\label{Lemma I} If $R$ is a $Gpf$-ring, then $R_{\mathfrak{p}}$ is a primary ring for all $\mathfrak{p}\in \Spec(R)$.\\
\end{lemma}

{\bf Proof.} Assume $(a/s)(b/t)=0$ and $a/s\neq0$ in $R_{\mathfrak{p}}$. Then there exists $u\in R\setminus\mathfrak{p}$ such that $uab=0$. But $\Ann(u^{n}b^{n})$ is a pure ideal for some $n\geqslant1$. Then by Lemma \ref{Remark V},  $\Ann(a) + \Ann(u^{n}b^{n})=R$. Thus $\Ann(u^{n}b^{n})\nsubseteq\mathfrak{p}$, since $\Ann(a)\subseteq\mathfrak{p}$. Therefore $b^{n}/t^{n}=0$ and so $R_{\mathfrak{p}}$ is a primary ring. $\Box$\\

\begin{proposition}\label{Proposition V} Every $Gpf$-ring is a mid ring.\\
\end{proposition}

{\bf Proof.} It follows from Theorem \ref{Theorem VII} and Lemma \ref{Lemma I}. $\Box$ \\

By $Q(R)$ we mean the totally ring fractions of $R$.\\

\begin{theorem}\label{Theorem IX} Let $R$ be a ring. Then the following are equivalent.\\
$\mathbf{(i)}$ $R$ is a p.p. ring.\\
$\mathbf{(ii)}$ $R$ is a mid ring and $Q(R)$ is a von Neumann regular ring.\\
\end{theorem}

{\bf Proof.} $\mathbf{(i)} \Rightarrow \mathbf{(ii)}$ : Let $a \in R$. Then $\Ann(a)$ is projective and so there exists idempotent $e \in R$ such that $\Ann(a)=Re$. Now, $a+e$ is a non-zero divisor. Then $Q(R)$ is a von Neumann regular ring.\\
$\mathbf{(ii)} \Rightarrow \mathbf{(i)}$ : Let $a \in R$. Then $\Ann(a)$ is a N-pure ideal. Since $Q(R)$ is regular, then there exists a non-zero divisor $s \in R$ such that $sa=a^{2}$. Thus there exist $n\geqslant1$ and $e \in \Ann(a)$ such that $(s-a)^{n}(1-e)=0$. Hence, we have $s^{n}(1-e)=ac$ for some $c \in R$. Then $e$ is an idempotent. But we have
$\Ann(a)=Re$ and so $R$ is a p.p. ring. $\Box$ \\

\begin{corollary}\label{Corollary V} Let $R$ be a ring. Then the following are equivalent.\\
$\mathbf{(i)}$ $R$ is a p.p. ring.\\
$\mathbf{(ii)}$ $R$ is a $Gpf$-ring and $Q(R)$ is a von Neumann regular ring.\\
\end{corollary}

{\bf Proof.} It follows from Proposition \ref{Proposition V} and Theorem \ref{Theorem IX}. $\Box$ \\

\begin{remark}\label{Remark I} Clearly every zero-dimensional ring is a Gpp-ring. Let $R$ be a zero-dimensional ring. Then for each $a\in R$ there exist $n\geqslant1$ and $b\in R$ such that $a^{n}(1-a^{n}b)=0$. Then $\Ann(a^{n})=R(1-a^{n}b)$ where $1-a^{n}b$ is an idempotent element of $R$. Thus $\mathbb{Z}/n\mathbb{Z}$ is a Gpp-ring and also is a mid ring for each positive integer $n$. Obviously, every p.f. ring is a mid ring. But the converse is not necessarily true. If $n=p_{1}^{k_{1}}...p_{t}^{k_{t}}$ where some $k_{j}\geqslant1$, then $\mathbb{Z}/n\mathbb{Z}$ is a mid ring by Proposition \ref{Proposition V} which is not p.f. ring.\\
\end{remark}

\begin{theorem}\label{Theorem IV} Let $R$ be a mid ring and $\mathfrak{p}$ be a prime ideal of $R$. Then $\mathfrak{p}$ is N-pure if and only if  $\mathfrak{p}\in \Min(R)$.\\
\end{theorem}

{\bf Proof.} Let $\mathfrak{p}$ be a N-pure ideal of $R$. If $\mathfrak{q}$ is a prime ideal of $R$ such that strictly contained in $\mathfrak{p}$, then there exists $a\in \mathfrak{p} \setminus \mathfrak{q}$. Hence there are $b\in \mathfrak{p}$ and $n\geqslant1$ such that $a^{n}(1-b)=0$. Thus $1 \in \mathfrak{p}$ which is a contradiction. Therefore, $\mathfrak{p}$ is a minimal prime ideal of $R$. Conversely, let $\mathfrak{p}$ be a minimal prime ideal of $R$ and $a \in \mathfrak{p}$. Using Theorem \ref{Theorem I}, it suffices to show that $\Ann(a^{m}) + \mathfrak{p}= R$ for some $m\geqslant1$. Assume that $\Ann(a^{n}) + \mathfrak{p} \neq R$ for all $n\geqslant1$. Setting $I:=\sum\limits_{n\geqslant1} \big(\Ann(a^{n}) + \mathfrak{p}\big)$. Then there exists a maximal ideal $\mathfrak{m}$ of $R$ such that $I\subseteq\mathfrak{m}$. So $\Ann(a^{n})\subseteq\mathfrak{m}$ for all $n\geqslant1$. Hence $a/1\notin \mathfrak{N}(R_{\mathfrak{m}})=\mathfrak{p}R_{\mathfrak{p}}$ and so $a\notin \mathfrak{p}$ which is a contradiction. $\Box$\\

\begin{theorem}\label{Theorem X} Every mid ring is a mp-ring.\\
\end{theorem}

{\bf Proof.} It follows from Theorems \ref{Theorem II} and \ref{Theorem IV}. $\Box$\\

\begin{example}\label{Example I} The converse of the above theorem is not necessarily true. As a specific example, let $R$ be the polynomial ring $k[x,y,z]$ modulo $I=(x^{3}-yz)$ where $k$ is a field. If $\mathfrak{p}=(x,z)$, then we have $I\subseteq\mathfrak{p}$. But $\mathfrak{p}/I$ is a prime ideal of $R$, since $R/(\mathfrak{p}/I)\simeq k[y]$. Now we consider the ring $R_{\mathfrak{m}}$ where $\mathfrak{m}=(\overline{x},\overline{y},\overline{z})$. Then  $\mathfrak{q}=(\overline{x}/1,\overline{z}^{2}/1)$ is a non-primary ideal of $R_{\mathfrak{m}}$ where $\overline{x}=x+I$. Because $(\overline{y}/1)(\overline{z}/1)=\overline{x}^{3}/1\in\mathfrak{q}$ but $\overline{y}/1\notin\sqrt{\mathfrak{q}}$
and $\overline{z}/1\notin\mathfrak{q}$. Therefore $R_{\mathfrak{m}}/\mathfrak{q}$ is a local quasi-prime ring and hence mp-ring which is not a mid ring.\\
\end{example}

We can deduce from Remark \ref{Remark I} and Example \ref{Example I}, that the class of mid rings is strictly between the class of p.f. rings(reduced mp-rings) and the class of mp-rings. Indeed,
$$\mathbf{p.f.\hskip 0.1cm rings=reduced\hskip 0.15cm mp-rings \subsetneq mid \hskip 0.15cm rings \subsetneq mp-rings}$$.\\

\end{document}